\documentclass[12pt]{amsart}
\usepackage{amsmath}
\usepackage{amscd}
\usepackage[psamsfonts]{amssymb}
\usepackage[mathscr,psamsfonts]{eucal}
\DeclareMathAlphabet{\eurm}{U}{eur}{m}{n}
\newcommand{\M}[1]{{\mathscr{#1}}}

\newcommand{\f}[1]{{\boldsymbol{#1}}}
\DeclareMathOperator{\id}{{id}}

\DeclareMathOperator{\proj}{{pr}}
\DeclareMathOperator{\lin}{{lin}}

\DeclareMathOperator{\Alt}{{Alt}}
\DeclareMathOperator{\pol}{pol}
\DeclareMathOperator{\inv}{inv}
\DeclareMathOperator{\Ker}{Ker}
\DeclareMathOperator{\Cla}{Cla}
\DeclareMathOperator{\pr}{pr}
\DeclareMathOperator{\byd}{\,{\raisebox{.1ex}{$\eurm :$}{\eurm =}}\,}
\newcommand{\sig}{\sigma}\newcommand{\del}{\delta}
\newcommand{\alp}{\alpha}
\newcommand{\bet}{\beta}\newcommand{\eps}{\epsilon}
\newcommand{\lam}{\lambda}
\newcommand{\ome}{\omega}
\newcommand{\Lam}{\Lambda}

\newcommand{\gam}{\gamma}
\newcommand{\bEq}{\begin{eqnarray}}
\newcommand{\eEq}{\end{eqnarray}}
\newcommand{\beq}{\begin{eqnarray*}}
\newcommand{\eeq}{\end{eqnarray*}}

\newcommand{\ac}[1]{\acute{#1}}

\newcommand{\dt}[1]{{\dot{#1}}}
\newcommand{\uten}[1]{\underset{#1}{\otimes}}
\newcommand{\ucar}[1]{\underset{#1}{\times}}
\newcommand{\sep}[1]{\qquad\text{#1}\qquad}
\newcommand{\ba}[1]{{{\bar{#1}}}}
\newcommand{\col}[3]{_{#1}{}^{#2}{}_{#3}}
\newcommand{\der}{\partial}
\newcommand{\nab}{\nabla}

\newcommand{\com}{\circ}
\newcommand{\car}{\times}
\newcommand{\ten}{\otimes}
\newcommand{\wed}{\wedge}
\newcommand{\Rn}{{I\!\!R}}

\begin{document}

\newcounter{theorem}

\newtheorem{definition}[theorem]{Definition}
\newtheorem{lemma}[theorem]{Lemma}
\newtheorem{proposition}[theorem]{Proposition}
\newtheorem{theorem}[theorem]{Theorem}
\newtheorem{corollary}[theorem]{Corollary}
\newtheorem{remark}[theorem]{Remark}
\newtheorem{example}[theorem]{Example}
\newtheorem{Note}[theorem]{Note}
\newcounter{assump}
\newtheorem{Assumption}{\indent Assumption}[assump]
\renewcommand{\thetheorem}{\thesection.\arabic{theorem}}

\newcommand{\bCr}{\begin{corollary}}
\newcommand{\eCr}{\end{corollary}}
\newcommand{\bDf}{\begin{definition}\em}
\newcommand{\eDf}{\end{definition}}
\newcommand{\bLm}{\begin{lemma}}
\newcommand{\eLm}{\end{lemma}}
\newcommand{\bPr}{\begin{proposition}}
\newcommand{\ePr}{\end{proposition}}
\newcommand{\bRm}{\begin{remark}\em}
\newcommand{\eRm}{\end{remark}}
\newcommand{\bEx}{\begin{example}\em}
\newcommand{\eEx}{\end{example}}
\newcommand{\bTh}{\begin{theorem}}
\newcommand{\eTh}{\end{theorem}}
\newcommand{\bNt}{\begin{Note}\em}
\newcommand{\eNt}{\end{Note}}
\newcommand{\bPf}{\begin{proof}[\noindent\indent{\sc Proof}]}
\newcommand{\ePf}{\end{proof}}


\title[Higher order reduction theorems]{Higher order reduction
theorems
\\
for classical connections and natural
\\
(0,2)-tensor fields
on the cotangent bundle}

\author{Josef Jany\v ska}



\address{
\newline
{Department of Mathematics, Masaryk University
\newline
Jan\'a\v ckovo n\'am. 2a, 662 95 Brno, Czech Republic}
\newline
E-mail: {\tt janyska@math.muni.cz}
}
\thanks{This paper has been supported
by the Ministry of Education of the Czech Republic under the Project
MSM 143100009.}

\begin{abstract} We generalize reduction theorems for classical
connections to operators with values in  $k$-th order natural bundles.
Using the first reduction theorem in order two we classify
all (0,2)-tensor fields on the cotangent bundle of
a manifold with a linear (non-symmetric) connection.
\end{abstract}
\maketitle

\section{Introduction}
\setcounter{equation}{0}
\setcounter{theorem}{0}
It is well known that natural operators of linear symmetric
connections on manifolds and of tensor fields
which have values in bundles of geometrical objects of order one
can be factorized through the curvature tensors, the tensor
fields and their covariant differentials. These results are known as
the first (the operators of connections only) and the second reduction
theorems.  The history of the first reduction theorem comes back to
the paper by Christoffel, \cite{Chr69}, and the history of the second
reduction theorem comes back to the paper
by Ricci and Levi Civita, \cite{RicLev01}.
For further references see \cite{Lub72, Sch54, ThoMic27}.
In \cite{Sch54} the proof for algebraic operators
(concomitants) is given.
In \cite{KolMicSlo93} the first and the second reduction
theorems are proved for all natural differential operators
by using the modern approach of natural
bundles and natural differential operators,
\cite{KruJan90, Nij72, Ter78}.

In this paper we generalize the reduction theorems
for natural operators which have values in higher order
natural bundles.

As an example we discuss natural (0,2)-tensor fields on the cotangent
bundle of a manifold.

In this paper we use the terms "natural bundle" and
``natural operator"
in the  sense of \cite{KolMicSlo93, KruJan90, Nij72, Ter78}.
Namely, a natural operator is defined to be a system of local operators
$A_{\f M}: C^\infty (F\f M)\to C^\infty (G\f M)$, such that
$A_{\f N}(f^*_Fs) = f^*_GA_{\f M}(s)$ for any section
$(s:\f M\to F\f M)\in
C^\infty(F\f M)$ and any (local) diffeomorphism
$f:\f M\to \f N$, where $F, G$ are two natural bundles and
$f^*_F s = Ff\com s\com f^{-1}$. A
natural operator is said to be of order $r$ if, for all sections
$s,q\in C^\infty (F\f M)$ and every point $x\in \f M$,
the condition $j^r_xs=j^r_xq$ implies
$A_{\f M}s(x)=A_{\f M}q(x)$. Then we have the induced natural
transformation
${\M A}_{\f M}:J^rF\f M\to G\f M$ such that
$A_{\f M}(s)={\M A}_{\f M}(j^rs)$, for all
$s\in C^\infty (F\f M)$. The correspondence
between natural operators of order $r$ and the induced natural
transformations is bijective. In this paper we shall identify
natural operators with the corresponding natural transformations.

Any natural bundle $F$ of order $r$ is given by its standard fibre
$S_F$ which is a left $G^r_m$-manifold, where
$G^r_m=\inv J^r_0(\Rn^m,\Rn^m)_0$ is the
$r$-th order differential group. A
classification of natural operators between natural bundles is
equivalent to the classification of equivariant maps between standard
fibers. Very important tool in classifications of equivariant maps is
the {\em orbit reduction theorem}, \cite{KolMicSlo93, Kru82, KruJan90}.
Let $p: G\to H$ be a Lie group epimorphism with the kernel $K$,
$M$ be a left $G$-space, $N,Q$ be left $H$-spaces and $\pi: M\to Q$ be
a $p$-equivariant surjective submersion, i.e. $\pi(gx) = p(g) \pi(x)$
for all $x\in M$, $g\in G$. Having $p$, we can consider every  left
$H$-space $N$ as a left $G$-space by $gy=p(g)y$, $g\in G$, $y\in N$.

\bTh\label{Th1.1}
If each $\pi^{-1}(q)$, $q\in Q$ is a $K$-orbit in $M$, then there
is a bijection between the $G$-maps $f: M\to N$ and the $H$-maps
$\varphi : Q\to N$ given by $f=\varphi\com \pi$.
\eTh

\section{Preliminaries}
\setcounter{equation}{0}
\setcounter{theorem}{0}

Let $\f M$ be an $m$-dimensional manifold.
If $(x^\lam)$, $\lam=1,\dots,m$,
is a local coordinate chart, then the induced coordinate charts on
$T\f M$ and $T^*\f M$
will be denoted by
$(x^\lam, \dot x^\lam)$ and $(x^\lam, \dot x_\lam)$ and the induced
local bases of sections of
$T\f M$ and
$T^*\f M$
will be denoted by
$(\der_\lam)$ and
$(d^\lam)$, respectively.

\bDf\label{Df2.1}
We define a {\em classical connection\/} to be a connection
\beq
\Lam : T\f M \to T^*\f M \uten{T\f M} TT\f M
\eeq
of the vector bundle
$p_{\f M}:T\f M \to \f M$,
which is linear and torsion free.
\eDf

The coordinate expression of a classical connection
$\Lam$
is of the type
\beq
\Lam = d^\lam \ten \big(\der_\lam
+ \Lam\col\lam\mu\nu \, \dot x^\nu \, \dt\der_\mu \big) \,,
\sep{with}
\Lam\col \mu \lam \nu = \Lam\col \nu \lam \mu
\in C^\infty(\f M,\Rn)
 \,.
\eeq

Classical connections can be regarded as sections of a 2nd
order natural bundle
$\Cla\f M \to \f M$, \cite{KolMicSlo93}.
The standard fibre of the functor $\Cla$ will be denoted by
${Q}=\Rn^{m}\ten \odot^2\Rn^{m*}$, elements of $Q$ will be said to
be {\em formal classical connections\/}, the
induced coordinates on ${Q}$ will be said to be {\em formal Christoffel
symbols\/} and will be
denoted by
$(\Lam\col \mu \lam \nu)$.

The action $\alp:G^2_m \car {Q} \to {Q}$
of the group $G^2_m$ on ${Q}$
is given in coordinates by
\beq
(\Lam\col \mu \lam \nu) \com \alp=
        (a^\lam_\rho \,( \Lam\col \sig \rho \tau
        \tilde a^\sig_\mu \tilde a^\tau_\nu
         - \tilde a^\rho_{\mu\nu}))\,
\eeq
where $(a^\lam_\mu, a^\lam_{\mu\nu})$ are the coordinates on
$G^2_m$ and $\tilde{\ }$ denotes the inverse element.

\bNt\label{Nt2.2}
Let us note that the action $\alp$ gives in a natural way
the action
\beq
\alp^r:G^{r+2}_m \car T^r_m{Q} \to T^r_m{Q} \,
\eeq
given by the jet prolongation of the action
$\alp$.
\eNt

\bRm\label{Rm2.3}
Let us consider the group epimorphism
$\pi^{r+2}_{r+1}: G^{r+2}_m\to  G^{r+1}_m$
and its kernel $B^{r+2}_{r+1}=\Ker \pi^{r+2}_{r+1}$. We have
the induced coordinates $(a^\lam_{\mu_1\dots\mu_{r+2}})$
on $B^{r+2}_{r+1}$. Then the restriction $\bar \alp^r$
of the action $\alp^r$ to $B^{r+2}_{r+1}$
has the following coordinate
expression
\begin{align*}
(\Lam\col {\mu_1} \lam {\mu_2}, & \dots,
        \Lam\col {\mu_1} \lam {\mu_2,\mu_3\dots\mu_{r+2}})
        \com \bar\alp^r
\\
 & =
        (\Lam\col {\mu_1} \lam {\mu_2},\dots,
        \Lam\col {\mu_1} \lam {\mu_2,\mu_3\dots\mu_{r+1}},
        \Lam\col {\mu_1} \lam {\mu_2,\mu_3\dots\mu_{r+2}} -
        \tilde a^\lam_{\mu_1\dots\mu_{r+2}})\,,
\end{align*}
where $(\Lam\col {\mu_1} \lam {\mu_2},
\Lam\col {\mu_1} \lam {\mu_2,\mu_3},  \dots,
        \Lam\col {\mu_1} \lam {\mu_2,\mu_3\dots\mu_{r+2}})$ are
the induced jet coordinates on $T^r_m{Q}$.
\eRm

The curvature tensor of a classical connection is a section
$
R[\Lam] :\f M\to W\f M\byd T^*\f M\ten T\f M\ten \wed^2T^*\f M
$
with coordinate expression
\beq
R[\Lam] = R\col \nu\rho{\lam\mu}\, d^\nu \ten
        \der_\rho \ten d^\lam \wed d^\mu \,,
\eeq
where the coefficients are
\beq
R_\nu{}^\rho{}_{\lam\mu} = \der_\mu \Lam\col \lam \rho \nu
        - \der_\lam \Lam\col \mu \rho \nu + \Lam\col \mu\sig \nu
        \Lam\col\lam\rho\sig - \Lam\col \lam\sig \nu
        \Lam\col\mu\rho\sig \,.
\eeq

Let us note that
the curvature tensor is a natural operator
\beq
R[\Lam]:C^\infty(\Cla\f M) \to C^\infty(W\f M)
\eeq
which is of order one, i.e., we have the associated
$G^3_m$-equivariant mapping, called the {\em formal curvature map\/}
of classical connections,
\beq
\M R:T^1_m{Q} \to W
\byd S_{T^*\ten T\ten
        \wed^2T^*}
\eeq
with coordinate expression
\bEq\label{Ex2.1}
(w_\nu{}^\rho{}_{\lam\mu})\com\M R= (\Lam\col \lam\rho{\nu,\mu}
        - \Lam\col \mu \rho {\nu,\lam} + \Lam\col \mu\sig \nu
        \Lam\col\lam\rho\sig - \Lam\col \lam\sig \nu
        \Lam\col\mu\rho\sig) \,,
\eEq
where $(w_\nu{}^\rho{}_{\lam\mu})$ are the induced
coordinates on the standard fibre $W
=\Rn^{m*}\ten\Rn^m\ten\wed^2\Rn^{m*}$.

\smallskip
Let $V\f M$ be a first order natural vector bundle over
(i.e., $V\f M$ is some tensor bundle over
$\f M$). Let us put $V_r\f M = V\f M\ten\ten^rT^*\f M$,
$V^{(k,r)}\f M = V_k\f M\ucar{\f M} \dots \ucar{\f M}V_r\f M$,
$V^{(r)}\f M\byd V^{(0,r)}\f M$.
Let us denote by
$V=\Rn^n$ or $V_r$ or $V^{(k,r)}$ the standard fibres of $V\f M$ or
$V_r\f M$ or $V^{(k,r)}\f M$, respectively.

The $r$-th order covariant differential of sections of $V\f M$ with
respect to classical connections is a natural operator
\beq
\nab^r:J^{r-1}\Cla\f M\ucar{\f M} J^rV\f M \to V_r\f M\,.
\eeq
We shall denote by the same symbol its corresponding
$G^{r+1}_m$-equivariant mapping
\beq
\nab^r:T^{r-1}_m Q \car T^r_m V \to V_r\,.
\eeq
We shall put
\beq\nab^{(k,r)}\byd (\nab^k,\dots,\nab^r):
J^{r-1}\Cla\f M\ucar{\f M} J^rV\f M \to V^{(k,r)}\f M\,
\eeq
and the same for the corresponding $G^{r+1}_m$-equivariant mapping.
Especially
$\nab^{(r)}\byd \nab^{(0,r)}$.

\bRm\label{Rm2.4}
For any section $\sig:\f M\to V\f M$ we have
\bEq\label{Ex2.2}
\Alt(\nab^2\sig) =
\pol(R[\Lam],
\sig)\,,
\eEq
where $\Alt$ is the antisymmetrization and
$\pol(R[\Lam],\sig)$ is a bilinear polynomial.
Namely, $\Alt(\nab^2R[\Lam])$ is a quadratic polynomial of $R[\Lam]$.
\eRm

If $(v^A)$ are
coordinates on $V$,
then $(v^A,v^A{}_\lam,\dots,v^A{}_{\lam_1\dots\lam_r})$ are the induced
jet coordinates on $T^r_m V$ (symmetric in all subscripts) and
$(V^A{}_{\lam_1\dots\lam_r})$ are the canonical coordinates
on $V_r$, then $\nab^r$ is of the form
\begin{align}\label{covariant differential}
(& V^A{}_{\lam_1\dots\lam_r})\com\nab^r
=v^A{}_{\lam_1\dots\lam_r} +\pol(T^{r-1}_m {Q} \car
        T^{r-1}_m V)\,,
\end{align}
where $\pol$ is a quadratic homogeneous polynomial on
$T^{r-1}_m {Q}\car  T^{r-1}_m V$.

\bRm\label{Rm2.5}
Let us recall the 1st and the 2nd Bianchi identities of classical
connections given in coordinates by
\beq
R\col {(\nu}\rho{\lam\mu)} = 0\,,\qquad
R\col {\nu}\rho{(\lam\mu;\sig)} = 0\,,
\eeq
respectively, where $;$ denotes the covariant differential with respect
to
$\Lam$ and $(\dots)$  denotes the cyclic
permutation.
\eRm

\section{The first $k$-th order reduction theorem}
\setcounter{equation}{0}
\setcounter{theorem}{0}

Let us introduce the following notations.

Let $W_0\f M\byd W\f M$,
$W_i\f M= W\f M\ten\ten^iT^*\f M$, $i\ge 0$. Let us put
$W^{(k,r)}\f M=
W_k\f M\car_{\f M} \dots \car_{\f M} W_r\f M$, $k\le r$.
We put $W^{(r)} \f M = W^{(0,r)} \f M$.
Then $W_i\f M$ and $W^{(k,r)}\f M$ are natural bundles
of order one and
the corresponding standard
fibers will be denoted by $W_i$ and $W^{(k,r)}$, respectively,
where $W_0\byd W$,
$W_i= W\ten\ten^i \Rn^{m*}$, $i\ge 0$, and $W^{(k,r)}=
W_k\car\dots \car W_r$.

We denote by
\beq
\M R_{i} : T^{i+1}_m {Q}\to  W_i
\eeq
the $G^{i+3}_m$-equivariant map associated with the
$i$-th covariant differential of curvature tensors of classical
connections
\beq
\nab^i R[\Lam]:C^\infty(\Cla\f M)\to C^\infty(W_i\f M)\,.
\eeq
The map $\M R_{i}$ is said to be the {\em formal
curvature map of order $i$} of classical connections.

Let $C_{i}\subset W_i$ be a subset given by
identities of the $i$-th covariant differentials of
the curvature tensors of classical connections,
i.e., by covariant differentials of the Bianchi identities
and the antisymmetrization of second order covariant differentials,
see Remark \ref{Rm2.4} and Remark \ref{Rm2.5}.
So $C_i$ is given by the following system of equations
\begin{align}
& w\col {(\nu}{\rho}{\lam\mu)\sig_1\dots \sig_i} =0\,,
\\
& w\col {\nu}{\rho}{(\lam\mu\sig_1)\sig_2\dots \sig_i} =0\,,
\\
& w\col {\nu}{\rho}{\lam\mu\sig_1\dots[\sig_{j-1}
        \sig_j]\dots \sig_i}
        + \pol(W^{(i-2)})=0\,,
\end{align}
where $j=2,\dots, i$ and $[..]$ denotes the antisymmetrization.

Let us put $C^{(r)}=
C_{0}\car \dots \car C_{r}$
and denote by
$C^{(k,r)}_{r^{(k-1)}}$, $k\le r$, the fiber in
$r^{(k-1)}\in C^{(k-1)}$ of the canonical
projection $\pr^r_{k-1}: C^{(r)}\to C^{(k-1)}$. For $r < k$
we put $C^{(k,r)}_{r^{(k-1)}}=\emptyset$.
Let us note that there is an affine structure on
the projection
$\pr^r_{r-1}: C^{(r)}\to C^{(r-1)}$, \cite{KolMicSlo93}.

Then we put
\bEq\label{e5.1}
\M R^{(k,r)}\byd (\M R_{k},\dots,\M R_{r}):
T^{r+1}_m {Q}\to  W^{(k,r)}\,,\qquad
\M R^{(r)}\byd \M R^{(0,r)}\,,
\eEq
which has values,
for any $j^{r+1}_0 \gam\in T^{r+1}_m Q$,
in $C^{(k,r)}_{\M R^{(k-1)}(j^{k}_0\gam)}$.
In \cite{KolMicSlo93} it was proved that $C^{(r)}$
is a submanifold
in $W^{(r)}$ and the restriction of $\M R^{(r)}$
to $C^{(r)}$ is a surjective submersion.
Then we can consider the fiber product
$T^{k}_m Q\car_{C^{(k-1)}} C^{(r)}$ and denote it
by $T^{k}_m Q\car C^{(k,r)}$.

\smallskip

First we shall prove the technical

\bLm\label{Lm3.1}
If  $r+1\ge k\ge 0$, then
the restricted
map
\beq
(\pi^{r+1}_{k},\M R^{(k,r)}):
T^{r+1}_m {Q}\to T^{k}_m Q\car C^{(k,r)}
\eeq
is a surjective submersion.
\eLm

\bPf
To prove surjectivity of $(\pi^{r+1}_{k},\M R^{(k,r)})$
it is sufficient to consider the commutative diagram
\beq
\begin{CD}
T^{r+1}_m Q @>\M R^{(r)}>> C^{(r)}
\\
@V{\pi^{r+1}_{k}}VV @VV{\pr^r_{k-1}}V
\\
T^{k}_m Q @>\M R^{(k-1)}>> C^{(k-1)}
\end{CD}
\eeq
All morphisms in the above diagram are surjective submersions
which implies that for any element $j^{k}_0\gam \in T^{k}_m Q$
the restriction of $\M R^{(r)}$ to the fibre
$(\pi^{r+1}_{k})^{-1}(j^{k}_0\gam)$ is a surjective submersion
of the fibre
$(\pi^{r+1}_{k})^{-1}(j^{k}_0 \gam)$
on the fibre  $(\pr^r_{k-1})^{-1}(\M R^{(k-1)}(j^k_0\gam))
\equiv C^{(k,r)}_{\M R^{(k-1)}(j^k_0\gam)}$
which proves that the mapping $(\pi^{r+1}_{k},\M R^{(k,r)})$
is surjective. To prove that $(\pi^{r+1}_{k},\M R^{(k,r)})$
is a submersion we shall consider the above diagram for $k=r$.
From the formal covariant differentials of (\ref{Ex2.1})
it follows, that $\M R^{(r,r)}= \M R_r$ is an affine
morphism over $\M R^{(r-1)}$ (with respect to the
affine structures on
$\pi^{r+1}_r:T^{r+1}_m Q\to T^r_m Q$ and
$\pr^r_{r-1}:C^{(r)}\to C^{(r-1)}$) which has a constant rank.
So the surjective morphism
$
(\pi^{r+1}_r, \M R_r): T^{r+1}_m Q\to T^r_m Q\car C^{(r,r)}
$
has a constant rank and hence is a submersion.
$
(\pi^{r+1}_{k},\M R^{(k,r)})
$
is then a composition of surjective submersions.
\ePf

Let $F$ be a natural bundle of order $k\ge 1$, i.e.,
$S_F$ is a left $G^k_m$-manifold.

\bTh\label{Th3.2}
Let $r+2\ge k$.
For every $G^{r+2}_m$-equivariant map
$$
f:T^r_m Q \to S_F
$$
there exists a unique
$G^k_m$-equivariant map $g:T^{k-2}_m Q \car C^{(k-2,r-1)}
\to S_F$ satisfying
\beq
f=g\com (\pi^r_{k-2}, \M R^{(k-2,r-1)})\,.
\eeq
\eTh

\bPf
Let us consider the space
$
S_{r} \byd \Rn^m\ten \odot^r\Rn^{m*}\,
$
with coordinates $(s^\lam{}_{\mu_1\mu_2\dots\mu_{r}})$.
Let us consider
the action of $G^r_m$ on $S_{r}$ given by
\bEq\label{e5.7}
\bar s^\lam{}_{\mu_1\mu_2\dots\mu_{r}}=
         s^\lam{}_{\mu_1\mu_2\dots\mu_{r}}
        - \tilde a^\lam_{\mu_1\dots\mu_{r}}
\,.
\eEq
From Remark \ref{Rm2.3} and (\ref{e5.7}) it is easy to see that
the symmetrization map $\sig_{s}:T^{r}_m {Q} \to  S_{r+2}$
given by
\beq
(s^\lam{}_{\mu_1\mu_2\dots\mu_{r+2}})\com
\sig_{s} = \Lam\col {(\mu_1}\lam{\mu_2,\mu_3\dots\mu_{r+2})}\,,
\eeq
is equivariant.

We have the $G^{r+2}_m$-equivariant
map
\beq
\varphi_{r} & \byd(\sig_{r},\pi^{r}_{r-1},
        \M R_{r-1})
     :  T^{r}_m {Q}  \to
        S_{r+2}  \car  T^{r-1}_m {Q}\car  W_{r-1}\,.
\eeq
On the other hand we define a $G^{r+2}_m$-equivariant map
\begin{align*}
\psi_{r} & : S_{r+2} \car
        T^{r-1}_m {Q}\car
        W_{r-1}\to
        T^{r}_m {Q}
        \,
\end{align*}
over the identity of $T^{r-1}_m {Q}$
by the following coordinate expression
\begin{align}
\Lam\col \mu\lam{\nu,\rho_1\dots\rho_{r}} & =
        s^\lam{}_{\mu\nu\rho_1\dots\rho_{r}}+\lin(
        w_\mu{}^\lam{}_{\nu\rho_1\dots\rho_{r}}
        -\pol(T^{r-1}_m Q))\,,\label{3}
\end{align}
where $\lin$ denotes a linear combination with real coefficients
which arises in the following way.
We recall that  $\M R_{r-1}$ gives the coordinate
expression, given by formal covariant differentials of (\ref{Ex2.1}),
\begin{align}
\Lam\col \mu\lam{\nu,\rho_1\dots\rho_{r}} -
        \Lam\col \mu\lam{\rho_1,\nu\rho_2\dots\rho_{r}}
        & = w_\mu{}^\lam{}_{\nu\rho_1\dots\rho_{s}} -
        \pol(T^{r-1}_m {Q})\,.\label{1}
\end{align}
We can write
\begin{align*}
\Lam\col \mu\lam{\nu,\rho_1\dots\rho_{r}} & =
        s^\lam{}_{\mu\nu\rho_1\dots\rho_{r}}+
        (\Lam\col \mu\lam{\nu,\rho_1\dots\rho_{r}} -
        \Lam\col {(\mu}\lam{\nu,\rho_1\dots\rho_{r})})\,.
\end{align*}
Then the term in brackets can be written as a linear
combination of terms of the type
\beq
\Lam\col \mu\lam{\nu,\rho_i\rho_1\dots\rho_{i-1}\rho_{i+1}\dots
        \rho_{r}} -
        \Lam\col \mu\lam{\rho_i,\nu\rho_1\dots
        \rho_{i-1}\rho_{i+1}\dots
        \rho_{r}}\,,
\eeq
$i=1,\dots,r$, and from (\ref{1}) we get (\ref{3}).

Moreover,
\beq
\psi_{r}\com\varphi_{r}=\id_{T^{r}_m {Q}}\,.
\eeq

Then the map $f\com \psi_{r}:
S_{r+2}\car
        T^{r-1}_m {Q}\car
        W_{r-1}\to S_F$ satisfies
the conditions of
the orbit reduction Theorem \ref{Th1.1} for the group epimorphism
$\pi^{r+2}_{r+1}:
G^{r+2}_m \to G^{r+1}_m$ and the surjective submersion
$\proj_{2,3}:S_{r+1} \car
        T^{r-1}_m {Q}\car
        W_{r-1}\to
        T^{r-1}_m {Q}\car
        W_{r-1}$.
Indeed, the space
$S_{r+2}$ is a $B^{r+2}_{r+1}$-orbit. Moreover,
(\ref{e5.7}) implies that the action
of $B^{r+2}_{r+1} $ on
$S_{r+2}$ is simply transitive.
Hence there exists a unique $G^{r+1}_m$-equivariant
map $g_{r-1}: T^{r-1}_m {Q}\car W_{r-1}
        \to S_F$
such that the following diagram
\begin{equation*}\begin{CD}
S_{r+2}  \car  T^{r-1}_m {Q}\car W_{r-1}
        @>{\psi_{r}}>> T^{r}_m {Q}
        @>f>>S_F
\\
@V{\proj_{2,3}}VV  @V{(\pi^{r}_{r-1},\M R_{r-1}
        )}VV @V\id_{S_F}VV
\\
        T^{r-1}_m {Q}\car
        W_{r-1} @>{\id_{T^{r-1}_m {Q}\car W_{r-1}}}>>
        T^{r-1}_m {Q}\car  W_{r-1}
        @>g_{r-1}>> S_F
\end{CD}\end{equation*}
commutes. So
$
f\com\psi_{r} = g_{r-1}\com \proj_{2,3}
$
and if we compose both sides with
$\varphi_{r}$, by
considering
$\proj_{2,3}\com\varphi_{r}
        = (\pi^{r}_{r-1},\M R_{r-1})$, we obtain
$f=g_{r-1} \com (\pi^{r}_{r-1},\M R_{r-1})\,.$

In the second step we consider the same construction for the
map $g_{r-1}$ and obtain the commutative diagram
\begin{equation*}\begin{CD}
(S_{r+1}  \car  T^{r-2}_m {Q}\car W_{r-2})\car W_{r-1}
        @>{\psi_{r-1}\car\id_{W_{r-1}}}>> T^{r-1}_m {Q} \car W_{r-1}
        @>g_{r-1}>>S_F
\\
@V{\proj_{2,3}\car\id_{W_{r-1}}}VV  @V{(\pi^{r-1}_{r-2},\M R_{r-2}
        )\car\id_{W_{r-1}}}VV @V\id_{S_F}VV
\\
        T^{r-2}_m {Q}\car W^{(r-2,r-1)}
        @>{\id_{T^{r-2}_m {Q}\car W^{(r-2,r-1)}}}>>
        T^{r-2}_m {Q}\car  W^{(r-2,r-1)}
        @>g_{r-2}>> S_F
\end{CD}\end{equation*}
So that there exists a unique $G^{r}_m $-equivariant
map
$
g_{r-2}: T^{r-2}\car W^{(r-2,r-1)}
        \to S_F
$
such that
$g_{r-1}= g_{r-2}\com
        ((\pi^{r-1}_{r-2},\M R_{r-2})\car\id_{W_{r-1}})
$,
i.e., $f=g_{r-2}\com (\pi^{r}_{r-2},\M R^{(r-2,r-1)})\,.$

Proceeding in this way we get in the last step
a unique $G^{k}_m $-equivariant
map $g_{k-2}: T^{k-2}_m Q\car W^{(k-2,r-1)}
\to S_F$ such that
\beq
f=g_{k-2}\com (\pi^{r}_{k-2},\M R^{(k-2,r-1)})\,.
\eeq
Putting $g$ the restriction of $g_{k-2}$ to $T^{k-2}_m Q\car
C^{(k-2,r-1)}$ we prove Theorem \ref{Th3.2}.{}
\ePf

In the above Theorem \ref{Th3.2} we have find
a map $g$ which factorizes $f$, but we did not prove, that
$(\pi^{r}_{k-2},\M R^{(k-2,r-1)}):T^{r}_m {Q}\to
T^{k-2}_m Q\car C^{(k-2,r-1)}$ satisfy the orbit conditions,
namely we did not prove that
$(\pi^{r}_{k-2},\M R^{(k-2,r-1)})^{-1}(j^{k-2}_0\lam,r^{(k-2,r-1)})$ is
a $B^{r+2}_{k}$-orbit for any $(j^{k-2}_0\lam,r^{(k-2,r-1)})
\in T^{k-2}_m Q\car C^{(k-2,r-1)}$. Now we shall
prove it.

\bLm\label{Lm3.3}
If $(j^r_0 \gam),(j^r_0 \ac{\gam}) \in
T^r_m {Q} $
satisfy
\beq
(\pi^r_{k-2},\M R^{(k-2,r-1)})(j^r_0 \gam)=
(\pi^r_{k-2},\M R^{(k-2,r-1)})(j^r_0\ac{\gam})\,,
\eeq
then there is an element
$h\in B^{r+2}_{k}$ such that
$h\,.\, (j^r_0 \ac{\gam}) = (j^r_0 \gam)$.
\eLm

\bPf
Consider the orbit set
$T^r_m {Q}/B^{r+2}_{k}$.
This is a $G^k_m $-set. Clearly the factor projection
\beq
p: T^r_m {Q} \to T^r_m {Q}/ B^{r+2}_{k}
\eeq
is a $G^{r+2}_m $-map. By Theorem
\ref{Th3.2} there is a $G^{k}_m$-equivariant
map
\beq
g: T^{k-2}_m Q\car C^{(k-2,r-1)}   \to
T^s_m {Q} /B^{r+2}_{k}
\eeq
satisfying $p=g\com (\pi^r_{k-2},\M R^{(k-2,r-1)})$.
If $(\pi^{r}_{k-2},\M R^{(k-2,r-1)})(j^r_0 \gam)=
(\pi^{r}_{k-2},\M R^{(k-2,r-1)})(j^s_0\ac{\gam})
=(j^{k-2}_0\lam, r^{(k-2,r-1)})$, then $p(j^r_0 \gam)
=g(j^{k-2}_0\lam,r^{(k-2,r-1)})
=p(j^r_0 \ac{\gam})$, so $j^r_0\gam$ and $j^r_0\ac{\gam}$ are
in the same $B^{r+2}_k$-orbit.
\ePf

It is easy to see that $T^{k-2}_m Q\car C^{(k-2,r-2)}$
is closed with respect to the action of the
group $G^k_m$. The corresponding
natural bundle of order $k$ is $J^{k-2}\Cla\f M\ucar{\f M}
C^{(k,r)}\f M$.  Then, as a
direct consequence of Theorem \ref{Th3.2}, we obtain the {\em first
$k$-order reduction theorem for classical connections}.

\bTh\label{Th3.4}
Let $F$ be a natural bundle of order $k\ge 1$
and let $r+2\ge k$. All
natural differential operators
$
f:C^\infty (\Cla\f M) \to  C^\infty(F\f E)
$
which are of order $r$ are of the
form
$$
f(j^r\Lam) = g(j^{k-2}\Lam,\nab^{(k-2,r-1)} R[\Lam])
$$
where $g$ is a unique natural operator
$$
g:J^{k-2}\Cla\f M\ucar{\f M} C^{(k-2,r-1)}\f M
        \to  F\f M\,.
$$
\eTh

\bRm\label{Rm3.5}
From the proof of Theorem \ref{Th3.2} it follows that the operator $g$
is the restriction of a natural operator defined on the
natural bundle $J^{k-2}\Cla\f M\ucar{\f M} W^{(k-2,r-1)}\f M$.
\eRm

\section{The second $k$-th order reduction theorem}
\setcounter{equation}{0}
\setcounter{theorem}{0}

(\ref{Ex2.2}) defines for $r=2$ the equation
$$
V^A_{[\lam\mu]} - \pol(C_0, V)=0\,\leqno{(E_2)}
$$
on $C_0\car V_2$ and for $r>2$ the system of equations
$$
V^A_{\mu_1\dots[\mu_{s-1} \mu_s]\dots \mu_r} -
\pol(C^{(r-2)},V^{(r-2)})=0\,
\leqno{(E_r)}
$$
on $C^{(r-2)}\car V^{(r)}$.

The $r$-th {\em Ricci subspace\/} $Z^{(r)}\subset
C^{(r-2)}\car V^{(r)}$ is defined by solutions of
$(E_2), \dots, (E_r)$, $r\ge 2$. For $r=0$ we put
$Z^{(0)}=V$ and for $r=1$ we put $Z^{(1)}= V^{(1)}$.
In \cite{KolMicSlo93} it was proved that $Z^{(r)}$
is a submanifold in $C^{(r-2)}\car V^{(r)}$ and
$(\M R^{(r-2)},\nab^{(r)})
   : T^{r-1}_m Q\car T^r_m V
\to  Z^{(r)}$ is a surjective submersion.
For $r > k-1$ we can consider the projection $\pr^r_{k-1}:
Z^{(r)}\to Z^{(k-1)}$ and denote by $Z^{(k,r)}_{z^{k-1}}$
its fiber in $z^{k-1}\in Z^{(k-1)}$.
Then we shall denote by
$
T^{k-2}_m Q\car T^{k-1}_m V \car Z^{(k,r)}
$
the fiber product
$
(T^{k-2}_m Q\car T^{k-1}_m V) \car_{Z^{(k-1)}} Z^{(r)}\,.
$

\bLm\label{Lm4.1}
If $r+1\ge k\ge 1$, then
the restricted map
\beq
(\pi^{r-1}_{k-2}\car\pi^r_{k-1})\car(\M R^{(k-2,r-2)},\nab^{(k,r)})
   : T^{r-1}_m Q\car T^r_m V
\to  T^{k-2}_m Q\car T^{k-1}_m V \car Z^{(k,r)}
\eeq
is a surjective submersion.
\eLm

\bPf
The proof of Lemma \ref{Lm4.1} follows from the commutative diagram
\beq
\begin{CD}
        T^{r-1}_m Q\car T^r_m V
        @>(\M R^{(r-2)},\nab^{(r)})>>
        Z^{(r)}
\\
        @V{\pi^{r-1}_{k-2}\car\pi^r_{k-1}}VV
        @VV{\pr^r_{k-1}}V
\\
        T^{k-2}_m Q \car T^{k-1}_m V
        @>(\M R^{(k-3)},\nab^{(k-1)})>>
        Z^{(k-1)}
\end{CD}
\eeq
where all morphisms are surjective submersions.
Hence
$
(\pi^{r-1}_{k-2}\car\pi^r_{k-1})\car(\M R^{(k-2,r-2)},\nab^{(k,r)})
$
is surjective. For $k=r$ the map $(\M R^{(r-2,r-2)}=\M R_{r-2},
\nab^{(r,r)}=\nab^r)$ is affine morphism over
$(\M R^{(r-3)}, \nab^{(r-1)})$ with a constant rank, i.e.,
$
(\pi^{r-1}_{r-2}\car\pi^r_{r-1})\car(\M R_{r-2},\nab^{r})
$
is a submersion.
$
(\pi^{r-1}_{k-2}\car\pi^r_{k-1})\car(\M R^{(k-2,r-2)},\nab^{(k,r)})
$
is then a composition of surjective submersions.
\ePf

\bTh\label{Th4.2}
Let $S_F$ be a left $G^k_m$-manifold. If $r+1\ge k\ge 1$, then
for every $G^{r+1}_m$-equivariant map $f: T^{r-1}_m Q\car T^r_m V
\to S_F$ there exists a unique $G^k_m$-equivariant map
$g: T^{k-2}_m Q\car T^{k-1}_m V \car Z^{(k,r)}\to S_F$
such that
$$
f = g\com
((\pi^{r-1}_{k-2}\car\pi^r_{k-1})\car(\M R^{(k-2,r-2)},\nab^{(k,r)}))\,.
$$
\eTh

\bPf
Consider the map
\beq
(\id_{T^{r-1}_m Q}\car \pi^{r}_{k-1},\nab^{(k,r)})
        : T^{r-1}_m {Q}\car T^{r}_m {V}
        \to T^{r-1}_m {Q}\car T^{k-1}_m {V} \car  V^{(k,r)}
\eeq
and denote by $\widetilde{V}^{(k,r)}\subset
        T^{r-1}_m {Q}\car T^{k-1}_m {V}
        \car V^{(k,r)}$
its image. By (\ref{covariant differential}), the restricted morphism
\beq
\tilde{\nab}^{(k,r)}:T^{r-1}_m {Q}\car T^{r}_m {V}
        \to
        \tilde{V}^{(k,r)}
\eeq
is bijective for every $j^{r-1}_0 \gam \in
T^{r-1}_m {Q}$, so that
$\tilde{\nab}^{(k,r)}$ is an equivariant diffeomorphism.
Define
\begin{align*}
\tilde{\M R}^{(k-2,r-2)} & :
\tilde{V}^{(k,r)}
 \to T^{k-2}_m Q\car T^{k-1}_m V\car Z^{(k,r)}
\end{align*}
by
\begin{align*}
\tilde{\M R}^{(k-2,r-2)}&
(j^{r-1}_0 \gam,j^{k-1}_0 {\mu},v) =
(j^{k-2}_0\gam, j^{k-1}_0 \mu, {\M
R}^{(k-2,r-2)}(j^{r-1}_0 \gam),v)\,,
\end{align*}
$(j^{r-1}_0 \gam,j^{k-1}_0 {\mu},v)\in
\tilde{V}^{(k,r)}$.
By Lemma \ref{Lm3.1}
$\tilde{\M R}^{(k-2,r-2)}$
is a surjective submersion.

Thus, Lemma \ref{Lm3.1} and Lemma \ref{Lm3.3} imply that
$\tilde{\M R}^{(k-2,r-2)}$ satisfies
the orbit conditions for the group epimorphism
$\pi^{r+1}_{k}: G^{r+1}_m \to G^{k}_m $ and
there exists a unique
$G^{k}_m$-equivariant map
$g: T^{k-2}_m Q\car T^{k-1}_m V \car  Z^{(k,r)}
\to S_F$
such that the diagram
\begin{equation*}\begin{CD}
\tilde{V}^{(k,r)}
        @>{(\tilde{\nab}^{(k,r)})^{-1}}>>
        T^{r-1}_m {Q}\car T^{r}_m {V}
        @>f>>
        S_F
\\
        @V{\tilde{\M R}^{(k-2,r-2)}}VV
        @V{(\pi^{r-1}_{k-2}\car \pi^{r}_{k-1},\nab^{(k,r)})}VV
        @V\id_{S_F}VV
\\
        T^{k-2}_m Q\car T^{k-1}_m V \car  Z^{(k,r)}
        @>{\id}>>
        T^{k-2}_m Q\car T^{k-1}_m V \car Z^{(k,r)}
        @>g>>
        S_F
\end{CD}\end{equation*}
commutes. Hence
$f\com (\tilde{\nab}^{(k,r)})^{-1}=g\com
\tilde{\M R}^{(k-2,r-2)}$.
Composing both sides with $\tilde{\nab}^{(k,r)}$,
by considering $\tilde{\M R}^{(k-2,r-2)}
\com \tilde{\nab}^{(k,r)}=(\pi^{r-1}_{k-2}\car \pi^{r}_{k-1},
\nab^{(k,r)})$,
we get
$$
f=g\com (\pi^{r-1}_{k-2}\car \pi^{r}_{k-1},\nab^{(k,r)})\,.\quad
$$
\vglue-\baselineskip
\ePf

$
T^{k-2}_m Q\car T^{k-1}_m V \car Z^{(k,r)}
$
is closed with respect to the action of the group $G^k_m$.
The corresponding natural bundle of order $k$ is
$
J^{k-2}\Cla\f M \ucar{\f M} J^{k-1}V\f M
        \ucar{\f M}Z^{(k,r)}\f M\,.
$

Then the second $k$-order reduction theorem can be
formulated as follows.

\bTh\label{Th4.3}
Let $F$ be a natural bundle of order $k\ge 1$ and let $r+1\ge k$. All
natural differential operators
$f:C^\infty(\Cla\f M\ucar{\f M} V\f M)\to C^\infty(F\f M)$
of order $r$ with respect sections of $V\f M$
are of the form
\beq
f(j^{r-1}\Lam,j^r\Phi) = g(j^{k-2}\Lam, j^{k-1}\Phi,\nab^{(k-2,r-2)}
R[\Lam],
        \nab^{(k,r)} \Phi)\,
\eeq
where $g$ is a unique natural operator
$$
g:J^{k-2}\Cla\f M\ucar{\f M} J^{k-1}V\f M\ucar{\f M}Z^{(k,r)}\f
M\to F\f M\,.
$$
\eTh

\bRm\label{Rm4.4}
The order $(r-1)$ of the above operators with respect to classical
connections is the minimal order we have to use. The second reduction
theorem can be easily generalized for any operators of order
$s\ge r-1$ with respect to connections. Then
$$
f(j^{s}\Lam,j^r\Phi) = g(j^{k-2}\Lam, j^{k-1}\Phi,\nab^{(k-2,s-1)}
R[\Lam],
        \nab^{(k,r)} \Phi)\,.
$$
\eRm

\bRm\label{Rm4.5}
If $\Lam$ is a linear non-symmetric connection on $\f M$, then there
exists its splitting $\Lam= \widetilde{\Lam} + T$, where
$\widetilde{\Lam}$ is the classical connection obtained by the
symmetrization of $\Lam$ and $T$ is the torsion tensor of $\Lam$.
Then all natural operators of order $r$ defined on $\Lam$ are of the
form
$$
f(j^r\Lam) = f(j^r \widetilde{\Lam}, j^rT) =
  g(j^{k-2}\widetilde{\Lam},j^{k-1}T,\widetilde{\nab}^{(k-2,r-1)}
R[\widetilde{\Lam}],
\widetilde{\nab}^{(k,r)}T)\,.
$$
\eRm

\bRm\label{Rm4.6}
If $g$ is a metric field on $\f M$, then there
exists the unique classical Levi Civita connection $\Lam$ given
by the metric field $g$.
Then, applying the second reduction theorem, we get that all natural
operators of order
$r\ge 1$ defined on
$g$ are of the form
$$
f(j^r g) = f(j^{r-1} {\Lam},j^r g) =
           h(j^{k-2}{\Lam},j^{k-1}g,\nab^{(k-2,r-2)}R[\Lam])
=h(j^{k-1} g,\nab^{(k-2,r-2)}R[\Lam]) \,.
$$
\eRm

\section{Natural (0,2)-tensor fields on the cotangent bundle}
\setcounter{equation}{0}
\setcounter{theorem}{0}

Typical applications of of higher order reduction
theorems are classifications of natural tensor fields on the tangent
(or  cotangent) bundle of a manifold endowed with a classical
connection or lifts of tensor fields to the tangent (or cotangent)
bundle by means of a classical connection,  see
\cite{Jan96, Jan01, KowSek88, Sek86, Sek88}.

As a direct consequence of Theorem \ref{Th3.2}, Theorem \ref{Th4.3}
and Remark \ref{Rm4.5} we get

\bCr\label{Cr5.1}
Let $(\f M,\Lam)$ be a manifold endowed with a linear (non-symmetric)
connection
$\Lam$. Then any natural tensor field $\Phi$ on $T\f M$
or $T^*\f M$ of order $r$
is of the type
\beq
\Phi(u,j^r \Lam) = \Phi(u,\widetilde\Lam, j^1T,
\widetilde\nab^{(r-1)} R[\widetilde\Lam],
\widetilde\nab^{(2,r)} T )\,,
\eeq
where $u\in T\f M$ or $u\in T^*\f M$, respectively, $\widetilde\Lam$
is the classical connection given by the symmetrization of $\Lam$ and
$T$ is the torsion tensor of $\Lam$.
\eCr


\bCr\label{Cr5.2}
Let $(\f M,\Lam, \Psi)$ be a manifold endowed with a linear
(non-symmetric) connection
$\Lam$ and a tensor field $\Psi$.
Then any natural tensor field
$\Phi$ on $T\f M$  or $T^*\f M$ of order $s$ with respect to $\Lam$
and of order $r$, $s\ge r-1$, with respect to $\Psi$ is of the type
\beq
\Phi(u,j^s \Lam, j^r\Psi) = \Phi(u,\widetilde\Lam, j^1 T, j^1\Psi,
\widetilde\nab^{(s-1)} R[\widetilde\Lam],
\widetilde\nab^{(2,s)} T,
\widetilde\nab^{(2,r)}\Psi)\,,
\eeq
where $u\in T\f M$ or $u\in T^*\f M$, respectively.
\eCr

As a concrete example let us classify all (0,2)-tensor fields
on $T^*\f M$ given by a linear (non-symmetric) connection $\Lam$.

\bTh\label{Th5.3}
Let $(\f M,\Lam)$ be a manifold endowed with a
linear (non-symmetric) connection
$\Lam$. Then all finite order natural (0,2)-tensor
fields on $T^*\f M$ are of order one and they form
a 14-parameter family of operators
with coordinate expression
\begin{align}\label{Ex5.0}
\Phi & = \big(A \,\dot x_\lam\, \dot x_\mu\,
+ C_1\, \dot x_\lam\,
T\col\rho\rho\mu + C_2\, \dot x_\mu\,
T\col\rho\rho\lam + C_3\, \dot x_\rho\,
T\col\lam\rho\mu
\\
&\quad
+ F_1\,
T\col{\rho}{\rho}{\lam}\,
T\col{\sig}{\sig}{\mu}
+ F_2\,
T\col{\sig}{\rho}{\lam}\,
T\col{\rho}{\sig}{\mu}
+ F_3\,
T\col{\rho}{\rho}{\sig}\,
T\col{\lam}{\sig}{\mu}\nonumber
\\
&\quad
+ G_1\,
T\col{\rho}{\rho}{\lam;\mu}
+ G_2\,
T\col{\rho}{\rho}{\mu;\lam}
+ G_3\,
T\col{\lam}{\rho}{\mu;\rho} \nonumber
\\
&\quad
+ H_1\, R\col \rho\rho{\lam\mu}
+ H_2\, R\col \lam\rho{\rho\mu}
\big)d^\lam\ten d^\mu\nonumber
\\
& \quad   +
 B \, d^\lam\ten (\dot{d}_\lam + \Lam\col \lam\rho\mu\, \dot x_\rho\,
 d^\mu) + C \, (\dot{d}_\lam + \Lam\col \lam\rho\mu\, \dot
x_\rho\,
 d^\mu) \ten d^\lam
\,,\nonumber
\end{align}
where $A, B, C, C_i, F_i, G_i, H_j$, $i=1,2,3$, $j=1,2$, are real constants.
\eTh

\bPf
Let us denote by
$S = \Rn^{m*}\car \ten^2 \Rn^{m*}\car \Rn^{m*}\ten \Rn^m\car
\Rn^m\ten \Rn^{m*} \car \ten^2\Rn^m $ the standard fibre
of
$\ten^2T^* (T^*\f M)$. The coordinates on $S$ will be denoted by
$(\dot x_\lam,{\phi}_{\lam\mu}, {\phi}_{\lam}{}^{\ba{\mu}} ,
{\phi}^{\ba{\lam}}{}_{{\mu}},
{\phi}^{\ba{\lam}\ba{\mu}})$. Then
we have the following action
of the group
$G^2_m$ on $S$
\begin{align*}
\ba{\dot{x}}_\lam & = \tilde a^\mu_\lam\, \dot{x}_\mu\,,
\\
\ba{\phi}_{\lam\mu} & = \tilde a^\rho_\lam\, \tilde a^\sig_\mu\,
\phi_{{\rho}{\sig}} +
a^\rho_\lam\, a^\alp_{\sig\bet}\, \tilde a^\bet_\mu \,
\tilde a^\kappa_\alp\, \dot{x}_\kappa\,\phi_{{\rho}}{}^{\ba{\sig}}
+
a^\sig_\mu\, a^\alp_{\rho\bet}\, \tilde a^\bet_\lam \,
\tilde a^\kappa_\alp\, \dot{x}_\kappa\,\phi^{\ba{\rho}}{}_{{\sig}}
\\
&\quad +
a^\alp_{\rho\bet}\, \tilde a^\bet_\lam \,
\tilde a^\kappa_\alp\, \dot{x}_\kappa\,
a^\gam_{\sig\del}\, \tilde a^\del_\mu \,
\tilde a^\nu_\gam\, \dot{x}_\nu\,
\phi^{\ba{\rho}\ba{\sig}}\,,
\\
\ba{\phi}_{\lam}{}^{\ba{\mu}} & =\tilde a^\rho_\lam\,
a^\mu_\sig\,
\phi_{{\rho}}{}^{\ba{\sig}} +
a^\mu_\sig\, a^\alp_{\rho\bet}\, \tilde a^\bet_\lam \,
\tilde a^\kappa_\alp\, \dot{x}_\kappa\,\phi^{\ba{\rho}\ba{\sig}}\,,
\\
\ba{\phi}^{\ba{\lam}}{}_{{\mu}} & = a^\lam_\rho\, \tilde a^\sig_\mu\,
\phi^{\ba{\rho}}{}_{\sig} +
a^\lam_\rho\, a^\alp_{\sig\bet}\, \tilde a^\bet_\mu \,
\tilde a^\kappa_\alp\, \dot{x}_\kappa\,\phi^{\ba{\rho}\ba{\sig}}\,,
\\
\ba{\phi}^{\ba{\lam}\ba{\mu}} & = a^\lam_\rho\, a^\mu_\sig\,
       \phi^{\ba{\rho}\ba{\sig}}\,.
\end{align*}

First let us discuss ${\phi}^{\ba{\lam}\ba{\mu}}$. We have,
by Corollary \ref{Cr5.1},
$$
{\phi}^{\ba{\lam}\ba{\mu}}= {\phi}^{\ba{\lam}\ba{\mu}}(\dot x_\lam,
 \widetilde\Lam\col \mu\lam\nu, T\col \mu\lam\nu, T\col \mu\lam{\nu,\sig},
 \widetilde{R}\col\nu\rho{\lam\mu;\sig_1;\dots;\sig_{i}},
 T\col\mu\lam{\nu;\sig_1;\dots;\sig_{j}})\,,
$$
$i=0,\dots, r-1$, $j=2,\dots, r$. The equivariance with
respect to homotheties $(c\,\del^\lam_\mu)$
implies
$$
c^2\,{\phi}^{\ba{\lam}\ba{\mu}}= {\phi}^{\ba{\lam}\ba{\mu}}(
c^{-1}\,\dot x_\lam,
c^{-1}\, \widetilde\Lam\col \mu\lam\nu,
c^{-1}\, T\col \mu\lam\nu,
c^{-2}\, T\col \mu\lam{\nu,\sig},
c^{-(i+2)}\, \widetilde R\col\nu\rho{\lam\mu;\sig_1;\dots;\sig_{i}},
c^{-(j+1)}\, T\col \mu\lam{\nu;\sig_1;\dots;\sig_j}
)
$$
which implies, by the homogeneous function theorem, \cite{KolMicSlo93},
that ${\phi}^{\ba{\lam}\ba{\mu}}$ is a polynomial of orders $a$ in
$\dot x_\lam$, $b$ in $\widetilde\Lam\col \mu\lam\nu$, $c_0$
in $T\col \mu\lam\nu$,  $c_1$
in $T\col \mu\lam{\nu,\sig}$, $d_i$ in
$\widetilde R\col\nu\rho{\lam\mu;\sig_1;\dots;\sig_i}$
and $e_j$ in $T\col \mu\lam{\nu;\sig_1;\dots;\sig_j}$
such that
\bEq\label{Ex5.1}
2 = -a -b -c_0-2\,c_1- \sum_{i=0}^{r-1} (i+2)\,d_i
- \sum_{j=2}^{r} (j+1)\,e_j\,.
\eEq
The equation (\ref{Ex5.1}) has no solution in natural numbers, so
we get by the homogeneous function theorem that
${\phi}^{\ba{\lam}\ba{\mu}}$ is independent of
all variables
and so it have to be
absolute invariant, hence
\beq
{\phi}^{\ba{\lam}\ba{\mu}}=0\,.
\eeq

For ${\phi}_{\lam}{}^{\ba{\mu}}$ and ${\phi}^{\ba{\lam}}{}_{{\mu}}$
we get from the equivariancy with respect to the homotheties
$(c\,\del^\lam_\mu)$ that they are polynomials of orders satisfying
\bEq\label{Ex5.2}
0 = -a -b -c_0-2\,c_1- \sum_{i=0}^{r-1} (i+2)\,d_i
- \sum_{j=2}^{r} (j+1)\,e_j\,.
\eEq
So also ${\phi}_{\lam}{}^{\ba{\mu}}$ and ${\phi}^{\ba{\lam}}{}_{{\mu}}$
are independent of
all variables and they have to be
absolute invariant, hence
\bEq
{\phi}_{\lam}{}^{\ba{\mu}}= B\,\del^\mu_\lam\,,\qquad
{\phi}^{\ba{\lam}}{}_{{\mu}}= C\,\del^\lam_\mu\,.
\eEq

Finally ${\phi}_{\lam\mu} $ has to be a polynomial of orders satisfying
\bEq\label{Ex5.4}
-2 = -a -b -c_0-2\,c_1- \sum_{i=0}^{r-1} (i+2)\,d_i
- \sum_{j=2}^{r} (j+1)\,e_j\,.
\eEq
There are 8 possible solutions of (\ref{Ex5.4}):

$a= 2$ and the others exponents vanish;

$a=1,\, b=1$ and the others exponents vanish;

$a=1,\, c_0=1$ and the others exponents vanish;

$b=2$ and the others exponents vanish;

$b=1,\, c_0=1$ and the others exponents vanish;

$c_0=2$ and the others exponents vanish;

$c_1=1$ and the others exponents vanish;

$d_0=1$ and the others exponents vanish.

It implies that the maximal order of the operator is one and
${\phi}_{\lam\mu} $  is of the form
\begin{align*}
{\phi}_{\lam\mu} & = A^{\rho\sig}_{\lam\mu}\, \dot x_\rho\, \dot x_\sig
+ B^{\rho\ome\tau}_{\lam\mu\kappa}\, \dot x_\rho\,
\widetilde\Lam\col\ome\kappa\tau
+ C^{\rho\ome\tau}_{\lam\mu\kappa}\, \dot x_\rho\,
T\col\ome\kappa\tau
\\
& \quad + D^{\ome_1\tau_1\ome_2\tau_2}_{\lam\mu\kappa_1\kappa_2}\,
\widetilde\Lam\col{\ome_1}{\kappa_1}{\tau_1}\,
\widetilde\Lam\col{\ome_2}{\kappa_2}{\tau_2}
+ E^{\ome_1\tau_1\ome_2\tau_2}_{\lam\mu\kappa_1\kappa_2}\,
\widetilde\Lam\col{\ome_1}{\kappa_1}{\tau_1}\,
T\col{\ome_2}{\kappa_2}{\tau_2}
\\
& \quad + F^{\ome_1\tau_1\ome_2\tau_2}_{\lam\mu\kappa_1\kappa_2}\,
T\col{\ome_1}{\kappa_1}{\tau_1}\,
T\col{\ome_2}{\kappa_2}{\tau_2}
+G^{\ome\tau\eps}_{\lam\mu\kappa}\,
T\col\ome\kappa{\tau,\eps}
+H^{\ome\tau\eps}_{\lam\mu\kappa}\,
\widetilde R\col\ome\kappa{\tau\eps}\,,
\end{align*}
where $A^{\rho\sig}_{\lam\mu}, \dots,
H^{\ome\tau\eps}_{\lam\mu\kappa}$ are absolute invariant tensors,
i.e.,
\begin{align*}
{\phi}_{\lam\mu} & = A\, \dot x_\lam\, \dot x_\mu
  + B_1\, \dot x_\lam\,
\widetilde\Lam\col\rho\rho\mu + B_2\, \dot x_\mu\,
\widetilde\Lam\col\rho\rho\lam + B_3\, \dot x_\rho\,
\widetilde\Lam\col\lam\rho\mu
\\
&\quad
  + C_1\, \dot x_\lam\,
T\col\rho\rho\mu + C_2\, \dot x_\mu\,
T\col\rho\rho\lam + C_3\, \dot x_\rho\,
T\col\lam\rho\mu
\\
&\quad
+ D_1\,
\widetilde\Lam\col{\rho}{\rho}{\lam}\,
\widetilde\Lam\col{\sig}{\sig}{\mu}
+ D_2\,
\widetilde\Lam\col{\sig}{\rho}{\lam}\,
\widetilde\Lam\col{\rho}{\sig}{\mu}
+ D_3\,
\widetilde\Lam\col{\rho}{\rho}{\sig}\,
\widetilde\Lam\col{\lam}{\sig}{\mu}
\\
&\quad
+ E_1\,
\widetilde\Lam\col{\rho}{\rho}{\lam}\,
T\col{\sig}{\sig}{\mu}
+ E_2\,
\widetilde\Lam\col{\sig}{\rho}{\lam}\,
T\col{\rho}{\sig}{\mu}
+ E_3\,
\widetilde\Lam\col{\rho}{\rho}{\sig}\,
T\col{\lam}{\sig}{\mu}
\\
&\quad
+ E_4\,
\widetilde\Lam\col{\rho}{\rho}{\mu}\,
T\col{\sig}{\sig}{\lam}
+ E_5\,
\widetilde\Lam\col{\sig}{\rho}{\mu}\,
T\col{\rho}{\sig}{\lam}
+ E_6\,
\widetilde\Lam\col{\lam}{\rho}{\mu}\,
T\col{\rho}{\sig}{\sig}
\\
&\quad
+ F_1\,
T\col{\rho}{\rho}{\lam}\,
T\col{\sig}{\sig}{\mu}
+ F_2\,
T\col{\sig}{\rho}{\lam}\,
T\col{\rho}{\sig}{\mu}
+ F_3\,
T\col{\rho}{\rho}{\sig}\,
T\col{\lam}{\sig}{\mu}
\\
&\quad
+ G_1\,
T\col{\rho}{\rho}{\lam,\mu}
+ G_2\,
T\col{\rho}{\rho}{\mu,\lam}
+ G_3\,
T\col{\lam}{\rho}{\mu,\rho}
\\
&\quad + H_1\,
\widetilde R\col\rho\rho{\lam\mu} +
H_2\,
\widetilde R\col\lam\rho{\rho\mu} \,.
\end{align*}

The equivariancy with respect to $(\del^\lam_\mu, a^\lam_{\mu\nu})$
implies $B_1=B_2 =0, B_3= B+C, D_i=0, E_1=E_4 =0,
E_2=G_3, E_3=-G_3, E_5=-G_3, E_6=-(G_1+G_2)  $ and the others
coefficients are arbitrary.
Then
\begin{align*}
{\phi}_{\lam\mu} & = A\, \dot x_\lam\, \dot x_\mu + (B+C)\,
\Lam\col\lam\rho\mu \, \dot x_\rho  + C_1\, \dot x_\lam\,
T\col\rho\rho\mu + C_2\, \dot x_\mu\,
T\col\rho\rho\lam + C_3\, \dot x_\rho\,
T\col\lam\rho\mu
\\
&\quad
+ F_1\,
T\col{\rho}{\rho}{\lam}\,
T\col{\sig}{\sig}{\mu}
+ F_2\,
T\col{\sig}{\rho}{\lam}\,
T\col{\rho}{\sig}{\mu}
+ F_3\,
T\col{\rho}{\rho}{\sig}\,
T\col{\lam}{\sig}{\mu}
\\
&\quad
+ G_1\,(
T\col{\rho}{\rho}{\lam,\mu}
        -\widetilde\Lam\col \lam\rho\mu\, T\col \rho\sig\sig)
+ G_2\,
(T\col{\rho}{\rho}{\mu,\lam}
       -\widetilde\Lam\col \mu\rho\lam\, T\col \rho\sig\sig)
\\
&\quad
+ G_3\,
(T\col{\lam}{\rho}{\mu,\rho}
       +\widetilde\Lam\col \sig\rho\lam\, T\col \rho\sig\mu
       -\widetilde\Lam\col \rho\rho\sig\, T\col \lam\sig\mu
       -\widetilde\Lam\col \sig\rho\mu\, T\col \rho\sig\lam)
\\
&\quad +H_1\,
\widetilde R\col\rho\rho{\lam\mu} +
H_2\,
\widetilde R\col\lam\rho{\rho\mu} \,,
\\
{\phi}_{\lam}{}^{\ba{\mu}}& = B\,\del_\lam^\mu\,, \qquad
{\phi}^{\ba{\lam}}{}_{{\mu}}  = C\,\del^\lam_\mu\,,
\qquad {\phi}^{\ba{\lam}\ba{\mu}}=0\,,
\end{align*}
which is the equivariant mapping corresponding to (\ref{Ex5.0}).
\ePf

\bRm
Let us note that the canonical symplectic form $\ome$ of $T^*\f M$
is a special case of (\ref{Ex5.0}). Namely, for $C=-B\ne0$ and the other
coefficients vanish we get just the scalar multiple of $\omega
= d^\lam\ten\dot d_\lam - \dot d_\lam\ten d^\lam$.
\eRm

The invariant description of the tensor fields (\ref{Ex5.0})
is the following. We have the canonical Liouville
1-form on $T^*\f M$ given in
coordinates by
\beq
\theta = \dot x_\lam\, d^\lam\,.
\eeq
The operator standing by $A$ is then $\theta\ten\theta$.

$\Lam$ gives 3-parameter family of (1,2) tensor fields on $\f M$,
\cite{KolMicSlo93}, given by
\bEq
S(\Lam) = C_1\, I_{T\f M}\ten \hat T
        + C_2\, \hat T\ten I_{T\f M}+ C_3\, T \,,
\eEq
where $\hat T$ is the contraction of the torsion tensor
and $I_{T\f M}: M\to T\f M\ten T^*\f M$ is the identity tensor.
Then the the evaluation $\langle S(\Lam), u\rangle$
gives three operators standing by $C_1, C_2, C_3$.

The connection $\Lam$ defines naturally the following 8 parameter
family of (0,2)-tensor fields on $\f M$, \cite{KolMicSlo93},
given by
\begin{align*}
G(\Lam) & = F_1\, C^{12}_{13} (T\ten T) +
        F_2\, C^{12}_{31} (T\ten T) +
        F_3\, C^{12}_{12} (T\ten T)
\\
&\quad
+
        G_1\, C^{1}_{1} \widetilde\nab T +
        G_2\, \overline{C^{1}_{1}  \widetilde\nab T} +
        G_3\, C^{1}_{3} \widetilde\nab  T
\\
&\quad
      +  H_1\, C^1_1 R[\widetilde\Lam]
        + H_2\, C^1_2 R[\widetilde\Lam]
        \,,
\end{align*}
where $C^{ij}_{kl}$ is the contraction with respect to indicated indices
and $\overline{C^{1}_{1}  \widetilde\nab T}$
denotes the conjugated tensor obtained
by the exchange of subindices.
The second 8-parameter subfamily of operators from (\ref{Ex5.0}) is then
given by the pullback of $G(\Lam)$ to $T^*\f M$.

The last two operators are given by the vertical projection
$$
\nu[\Lam^*]: T^*\f M \to T^*T^*\f M\ten VT^*\f M\,,\quad
\nu[\Lam^*]=(\dot{d}_\lam + \Lam \col \lam\rho\mu\,
        \dot{x}_\rho\, d^\mu) \ten \dot{\der}^\lam
$$
associated with the connection $\Lam^*$ dual to $\Lam$. Two possible
contractions
of $\nu[\Lam^*]\ten \ome$ give the
last two operators.

\bRm
    For a symmetric connection $\Lam$ the family (\ref{Ex5.0})
reduces to 5-parameter family.
\eRm



\end{document}